\documentclass[10pt]{article}

\usepackage{amsfonts,latexsym,amstext}

\newtheorem{theorem}{Theorem}[section]

\newtheorem{proposition}[theorem]{Proposition}
\newtheorem{definition}{Definition}[section]

\newtheorem{hypothesis}[theorem]{Hypothesis}

\makeatletter
\@addtoreset{equation}{section}

\makeatother

\def\sqr#1#2{{\vcenter{\vbox{\hrule height .#2pt \hbox{\vrule
 width .#2pt height#1pt \kern#1pt \vrule
width .#2pt} \hrule height .#2pt}}}}

\def\ds{\begin{displaystyle}}
\def\eds{\end{displaystyle}}
\def\dis{\displaystyle }
\def\<{\langle }
\def\>{\rangle }

\def\dim{\noindent \hbox{{\bf Proof.} }}

\def\R{\mathbb R}

\def\E{\mathbb E}
\def\P{\mathbb P}

\def\cale{{\cal E}}
\def\calf{{\cal F}}

\def\calh{{\cal H}}

\def\call{{\cal L}}
\def\calo{{\cal O}}

\title{Stochastic maximum principle for optimal control of SPDEs.
\thanks{Supported by 	
Marie Curie ITN
Call: FP7-PEOPLE-2007-1-1-ITN, no. 213841-2:
\em{Deterministic and Stochastic Controlled Systems and Applications }}}

\date{}
\author{Marco Fuhrman
\\
Dipartimento di Matematica,
Politecnico di Milano\\
Piazza Leonardo da Vinci 32, 20133 MILANO, Italy\\
e-mail: marco.fuhrman@polimi.it
\\ \\
Ying Hu\\
IRMAR, Universit\'e Rennes 1
\\ Campus de Beaulieu, 35042 RENNES Cedex, France
\\ e-mail: ying.hu@univ-rennes1.fr
\\ \\
Gianmario Tessitore
\\
Dipartimento di Matematica, Universit\`a di Milano-Bicocca\\
Via R. Cozzi 53, 20125 MILANO, Italy\\
e-mail: gianmario.tessitore@unimib.it
}

\begin{document}

\maketitle



\noindent {\bf Abstract.} 

 In this note, we give the stochastic maximum principle for optimal control of stochastic PDEs in the general case (when the control domain need not be convex and the diffusion coefficient can contain a control variable).

\bigskip

\section{Introduction} The problem of finding necessary optimality conditions for stochastic optimal control problems (generalizing in this way the Pontryagin
maximum principle to the stochastic case) has been solved in great generality, in the classical finite dimensional case, in the well known paper by S. Peng \cite{Peng}. The author allows the set of control actions to be non-convex and the diffusion coefficient to depend on the control; consequently he is led to introducing the equations for the second variation process and for its dual. As far as  infinite dimensional equations are concerned the cases in which the control domain is convex or diffusion does not depend on the control have been treated in \cite{Ben, HuPe}. On the contrary in the general case (when the control domain need not be convex and the diffusion coefficient can contain a control variable) existing results  are limited to abstract evolution equations under assumptions that are not satisfied by the  large majority of concrete stochastic PDEs (for instance the case of Nemitsky operators  on $L^p$ spaces is not covered, see  \cite{TaLi, LuZh}). Here we formulate a controlled parabolic stochastic PDE in a semi-abstract way and show how the specific regularity properties on the semigroup corresponding to the differential operator can be used to treat the case of Nemitsky-type coefficients. The key point is the explicit characterization of the second  dual process  in term of a suitable quadratic functional, see Definition \ref{def-P}.

\section{Formulation of the optimal control problem}

Let $\calo
\subset \R^n$ be a bounded open set with regular boundary. We consider the following controlled SPDE formulated in a partially abstract way  in
the state space $H=L^2(\calo)$
(norm $|\cdot|$,  scalar product $\<\cdot,\cdot\>$):
\begin{equation}\label{state-equation}
    \left\{
\begin{array}{lll}
dX_t(x)&=&\dis AX_t(x)\,dt +b(x,X_t(x),u_t)\,dt+
\sum_{j=1}^m\sigma_j(x,X_t(x),u_t)\,d\beta^j_t,\qquad  t\in [0,T], x\in\calo,
\\
X_0(x)&=&x_0(x),
\end{array}
\right.
\end{equation}
and the following cost functional:
$$
J(u)=\E \int_0^T\int_\calo l(x,X_t(x),u_t)\,dx\,dt
+\E \int_\calo h(x,X_T(x))\,dx.
$$
We work in the following setting.

\begin{hypothesis}\label{standing_ass}

\begin{enumerate}
\item $A$ is the realization of a partial differential operator
with appropriate boundary conditions. We assume that $A$ is  the infinitesimal generator of
  a strongly continuous semigroup $e^{tA}$, $t\ge 0$,
in $H$.
Moreover,
for every $p\in [2,\infty)$ and $t\in [0,T]$,
$
e^{tA}(L^p(\calo))\subset L^p(\calo)$ with
$\|e^{tA}f\|_{L^p(\calo)}\le C_{p,T} \|f\|_{L^p(\calo)}
$
for some constants $C_{p,T}$ independent of $t$ and $f$.
Finally the restriction of $e^{tA}$, $t\ge 0$, to $L^4(\calo)$ is an analytic semigroup
with domain of the infinitesimal generator compactly embedded in $L^4(\calo)$.

\item $(\beta^1_t,\ldots,\beta^m_t)$, $t\ge 0$ is a standard
$m$-dimensional Wiener process on a complete probability space $(\Omega,\cale,\mathbb{P})$ and we denote by $(\calf_t)_{t\ge 0}$
its natural (completed) filtration. All stochastic processes will be progressively measurable with respect to
$(\calf_t)_{t\ge 0}$.

\item $b,\sigma_j (j=1,...,m), l:\calo\times\R\times U\to \R$ and $h:\calo\times\R\to \R$  are measurable functions.
We assume that they are   continuous with respect to the third  variable (the control variable),
of class $C^2$ with respect to the second (the state variable), and bounded
 together with their first and second derivative with respect to the second variable.

 \item the set of admissible control actions is a separable metric
space $U$ and an admissible control $u$ is a  (progressive) process with values in $U$.
\end{enumerate}
\end{hypothesis}

Under the above conditions, for every control $u$ there exists a unique
mild solution, i.e. a continuous process in $H$ such that, $\P$-a.s.
$$
X_t=e^{tA}x_0+\int_0^te^{(t-s)A} b(\cdot,X_s(\cdot),u_s)\,ds+
\int_0^te^{(t-s)A}\sigma_j(\cdot,X_s(\cdot),u_s)\,d\beta^j_s,
\quad t\in [0,T].
$$

%
%
%
\section{Expansions of the solution and of the cost}
We assume that an optimal control $\bar{u}$ exists and denote by $\bar{X}$ the corresponding optimal state.
We introduce the spike variation: we fix an arbitrary
interval $[\bar{t},\bar{t}+\epsilon]\subset (0,T)$ 
and an arbitrary $U$-valued, $\calf_{\bar{t}}$-measurable random variable $v$
define the following perturbation of $\bar u$:
$u_t^\epsilon=vI_{[\bar{t},\bar{t}+\epsilon]}(t)+
\bar{u}_t I_{[\bar{t},\bar{t}+\epsilon]^c}(t)
$ and denote by $X^\epsilon$ the solution of the state equation (\ref{state-equation}) with control $u=u^\epsilon$.

 We introduce two linear equations corresponding to first and second expansion
 of $X^{\epsilon}$ with respect to $\epsilon$
 (both equations are understood in the mild sense). In the following, derivatives with respect to the
 state variable will be denoted $b',b'',\sigma',\sigma''$ and
$$
\delta b_t(x)=b(x,\bar{X}_t(x),
u_t^\epsilon)-b(x,\bar{X}_t(x),
\bar{u}_t),\qquad
\delta \sigma_{jt}(x)=\sigma_j(x,\bar{X}_t(x),
u_t^\epsilon)-
\sigma_j(x,\bar{X}_t(x),\bar{u}_t),
$$
$$
\delta b_t'(x)=b'(x,\bar{X}_t(x),
u_t^\epsilon)-b'(x,\bar{X}_t(x),
\bar{u}_t),\qquad
\delta \sigma_{jt}'(x)=\sigma_j'(x,\bar{X}_t(x),
u_t^\epsilon)-
\sigma_j'(x,\bar{X}_t(x),\bar{u}_t).
$$
Consider
\begin{equation}\label{eq-Y}
\left\{
\begin{array}{lll}
dY^\epsilon_t(x)&=&
\bigg[AY^\epsilon_t(x) +b'(x,\bar{X}_t(x),\bar{u}_t)\cdot
Y^\epsilon_t(x)\bigg]\,dt
+
\sigma_j'(x,\bar{X}_t(x),\bar{u}_t)\cdot Y^\epsilon_t(x)
\,d\beta^j_t
+\delta b_t(x)\,dt
+\delta \sigma_{jt}(x)\,d\beta^j_t
\\
Y^\epsilon_0(x)&=&0
\end{array}
\right.
\end{equation}
\begin{equation}\label{eq_Z}
\left\{
\begin{array}{lll}
dZ^\epsilon_t(x)& \!\!\! \!\!\!=\!\!\!&\!\!\!\!\!\!\!\!\!
\bigg[AZ^\epsilon_t(x) +b'(x,\bar{X}_t(x),\bar{u}_t)\cdot
Z^\epsilon_t(x)\bigg]\,dt
+ \sigma_j'(x,\bar{X}_t(x),\bar{u}_t)\cdot Z^\epsilon_t(x)
\,d\beta^j_t
\\&& \!\!\!\!\!\!\!\!\!\!\!\!\!\!\!+\bigg[
\frac{1}{2}b''(x,\bar{X}_t(x),\bar{u}_t)\cdot
Y^\epsilon_t(x)^{{2}}+
\delta b_t'(x)\cdot Y^\epsilon_t(x)\bigg]dt
+\bigg[
\frac{1}{2}\sigma_j''(x,\bar{X}_t(x),\bar{u}_t)\cdot
Y^\epsilon_t(x)^{{2}}
+\delta \sigma_{jt}'(x)\cdot Y^\epsilon_t(x)\bigg]d\beta^j_t
\\
Z^\epsilon_0(x)&=&0
\end{array}
\right.
\end{equation}
We notice that  to formulate the second equation in $H$ we need to show that the first admits solutions in $L^4(\calo)$.

$ $

\noindent The following proposition states existence and uniqueness of the solution
to the above equation in all spaces $L^p(\calo)$ together with the estimate
of their dependence with respect to $\epsilon$.
 The proof is technical but based on standard estimates  and we omit it.
\begin{proposition}\label{estimates on YZ} Equations (\ref{eq-Y}) and (\ref{eq_Z}) admit a unique continuous mild solution. Moreover
for all $p\geq 2$
$$
\sup_{t\in[0,T]}\left(\sqrt{\epsilon}^{-1}(\E  \|Y^\epsilon_t\|_{L^p(\calo)}^p)^{1/p}+
\epsilon^{-1}(\E  \|Z^\epsilon_t\|_{L^p(\calo)}^p)^{1/p}\right)
\le C_p,\qquad
\sup_{t\in[0,T]} \left(\E \|
X^\epsilon_t-
\bar{X}_t-Y^\epsilon_t-
Z^\epsilon_t\|_H^2\right)^{1/2}
= o(\epsilon).
$$
\end{proposition}
As far as the cost is concerned we set
$
\delta l_t(x)=l(x,\bar{X}_t(x),u_t^\epsilon)-l(x,\bar{X}_t(x),\bar{u}_t)
$
and prove that
\begin{proposition}\label{first-expansion-cost}
$$
 J(u^\epsilon)-J(\bar{u})=
\E \int_0^T\int_\calo \delta l_t(x)\,dx\,dt
+\Delta_1^\epsilon+
\Delta_2^\epsilon
+  o(\epsilon),
$$
where
$$
\begin{array}{lll}
\Delta_1^\epsilon&=&\dis
\E \int_0^T\int_\calo l'(x,\bar{X}_t(x),\bar{u}_t)
(Y_t^\epsilon(x)+Z_t^\epsilon(x)) \normalcolor\,dx\,dt
+\E \int_\calo h'(x,\bar{X}_T(x))
(Y_T^\epsilon(x)+Z_T^\epsilon(x))
\normalcolor\,dx,
\\
\Delta_2^\epsilon
&=&\dis
\frac{1}{2}\, \E \int_0^T\int_\calo l''(x,\bar{X}_t(x),\bar{u}_t)
\,Y_t^\epsilon(x)^2\normalcolor \,dx\,dt
+\frac{1}{2}\,\E \int_\calo h''(x,\bar{X}_T(x))\,
Y_T^\epsilon(x)^2\normalcolor\,dx.
\end{array}
$$
\end{proposition}
\section{The first and second adjoint processes}
The following proposition is special case of a result in \cite{HuPe}:
\begin{proposition}\label{def-pq} Let $A^*$ be the $L^2(\calo)$-adjoint operator of $A$. There exists a unique $m+1$-tuple of $L^2(\calo)$ processes
$(p,q_{j})$, with $p$ continuous and $
 \E  \sup_{t\in[0,T]}|
p_t|^2+
\E \int_0^T |
q_{jt}|^2\,dt<\infty
$,
that verify (in a mild sense) the backward stochastic differential equation:
$$
\left\{
\begin{array}{lll}
-dp_t(x)\!\!\!&=&\!\!\!-q_{jt}(x)\,d\beta^j_t+
\bigg[A^*p_t(x) +b'(x,\bar{X}_t(x),\bar{u}_t)\cdot
p_t(x)
+
\sigma_j'(x,\bar{X}_t(x),\bar{u}_t)\cdot q_{jt}(x)
+l'(x,\bar{X}_t(x),\bar{u}_t)\bigg]
\,dt
\\
p_T(x)\!\!\!&=&\!\!\!h'(x,\bar{X}_T(x)).
\end{array}
\right.
$$
\end{proposition}
The following proposition formally follows from
Proposition \ref{first-expansion-cost}
 computing the It\^o differentials
$
d\int_\calo Y^\epsilon_t(x)p_t(x)\,dx$ and
$
d\int_\calo Z^\epsilon_t(x)p_t(x)\,dx,
$
while the formal proof goes through Yosida approximations of $A$.
\begin{proposition} We have
\begin{equation}\label{penultima_stima}
 J(u^\epsilon)-J(u)=
\E \int_0^T\int_\calo \big[\delta l_t(x)+
\delta b_t(x)
p_t(x)
+
\delta \sigma_{jt}(x)
q_{jt}(x)\big]
\,dx\,dt
+\frac{1}{2}\,\Delta_3^\epsilon
+  o(\epsilon),
\end{equation}
where
\begin{equation}\label{delta_3}
\Delta_3^\epsilon
=
 \E \int_0^T \int_\calo \bar{H}_t(x)
\,Y_t^\epsilon(x)^2 dxdt
+\E \int_\calo \bar{h}(x)\,Y_T^\epsilon(x)^2\,dx,
\end{equation}
with
$$
\bar{H}_t(x)=
l''(x,\bar{X}_t(x),\bar{u}_t)
 +b''(x,\bar{X}_t(x),\bar{u}_t)p_t(x)+
\sigma_j''(x,\bar{X}_t(x),\bar{u}_t)
q_{jt}(x),\quad
\bar{h}(x)=h''(x,\bar{X}_T(x)).
$$
We notice that the multiplication by $\bar{H}_t(\cdot)$ is not a bounded operator in $H$.
\end{proposition}

\begin{definition}\label{def-P}
For fixed $t\in[0,T]$ and $f\in L^4(\calo)$, we consider the equation (understood as usual in mild form)
\begin{equation}\label{equat-for-Y}
    \left\{\begin{array}{lll}
    dY^{t,f}_s(x)\!\!\!\!&=&\!\!\!\!\dis
    AY^{t,f}_s(x)\,ds + b'(x,\bar X_s(x),\bar u_s)Y^{t,f}_s(x)\,ds
     +  \sigma_j'(x,\bar X_s(x),\bar u_s)Y^{t,f}_s(x)\,dW^j_s,
     \qquad s\in [t,T],
    \\
    Y^{t,f}_t(x) \!\!\!\!&=&\!\!\!\!f(x).
\end{array}\right.
\end{equation}
We denote ${ \call}$ the space of bounded linear operators ${ L^4(\calo)\to L^4(\calo)^*=L^{4/3}(\calo)}$ and define
a progressive process $(P_t)_{t\in [0,T]}$ with values
in $\call$ setting for $t\in [0,T]$, $f,g\in L^4(\calo)$,
$${
    \< P_t f,g\>= \E^{\calf_t}\int_t^T\int_{\calo}\bar H_s(x) Y_s^{t,f}(x) Y_s^{t,g}(x)\, dx\,ds
+
 \E^{\calf_t}\int_{\calo} \bar
   h (x) Y_T^{t,f}(x)    Y_T^{t,g}(x) \, dx \qquad \P-a.s.}
 $$
 \end{definition}
(by abuse of language by $\<\cdot,\cdot\>$ we also denote the duality between $L^4(\calo)$ and $L^{4/3}(\calo)$).
\newline\noindent
Exploiting the analyticity of the semigroup generated by $ A$ on $L^4(\calo)$ we prove the following proposition that is the key point for our final argument.
\begin{proposition} \label{prop_stime-su-P} We have
 $
    \sup_{t\in [0,T]}\E\|P_t\|_\call^2<\infty
$. Moreover
$
\E|\<P_{t+\epsilon}-P_t)f,g\>|\to 0,$ as  $\epsilon \rightarrow 0$,  $\forall f,g\in L^4(\calo)$. Finally for every $\eta\in (0,1/4)$ there exists
a constant $C_\eta $ such that  \begin{equation}\label{stimaindaetaperp}
| \<P_t(-A)^\eta f,(-A)^\eta g\>| \le
C_\eta \|f\|_4\|g\|_4 (T-t)^{-2\eta} \left[
\left(\int_t^T\E^{\calf_t}|\bar H_s|^2ds\right)^{1/2}
+
\left(\E^{\calf_t}|\bar h|^2 \right)^{1/2}
\right]
,\; \P-a.s.
 \end{equation}
where $D(-A)^\eta$ is the domain of the fractional power of $A$ in $L^4(\calo) $ and by $ \|\,\cdot\,\|_4$
we denote the norm in $L^4(\calo) $.
\end{proposition}

\section{The Maximum Principle}
For $u\in U$ and $X,p,q_1,\ldots,q_m\in L^2(\calo)$  denote
$$
\calh(u,X,p,q_1,\ldots,q_m)=
\int_\calo \bigg[ l(x,X(x),u)+
 b(x,X(x),u)
p(x)
+
\sigma_{j}(x,X(x),u)
q_{j}(x)\bigg]
\,dx
$$
\begin{theorem} Let $(\bar{X}_t,\bar{u}_t)$ be an optimal pair
and let $p,q_1,\ldots,q_m$ be defined as in  Proposition \ref{def-pq}
and $P$ be defined as in Definition \ref{def-P}.
Then the following inequality holds $\P$-a.s. for a.e. $t\in
[0,T]$ and for every $v\in U$:
$$
\begin{array}{l}
\calh(v,\bar{X}_t,p_t, q_{1t},\ldots,q_{mt})-
 \calh(\bar{u}_t,\bar{X}_t,p_t, q_{1t},\ldots,q_{mt})
 \\\dis
 \qquad
+\frac{1}{2}\,\<P_t[\sigma_j(\cdot,\bar{X}_t(\cdot),v)-
\sigma_j(\cdot,\bar{X}_t(\cdot),\bar{u}_t)],
\sigma_j(\cdot,\bar{X}_t(\cdot),v)-
\sigma_j(\cdot,\bar{X}_t(\cdot),\bar{u}_t)\>
\ge 0.
\end{array}
$$
\end{theorem}

\dim
By the Markov property of the solutions to equation (\ref{equat-for-Y}) and Proposition \ref{estimates on YZ}
we get:
\begin{equation}\label{eq-penultima}
 \begin{array}{rcl}\dis
\E\int_0^T \!\!\!\<\bar H_s  Y_s^\epsilon,Y_s^\epsilon \> \,ds
+
 \E\<
 \bar h Y_T^\epsilon,Y_T^\epsilon \>
 \!\!\!&=&\!\!\!
\dis \E\int_{t_0+\epsilon}^T \!\!\! \!\!\!\<\bar H_s
 Y_s^{t_0+\epsilon,Y_{t_0+\epsilon}^\epsilon}
,Y_s^{t_0+\epsilon,Y_{t_0+\epsilon}^\epsilon} \> \,ds
+
 \E\<
 \bar h Y_T^{t_0+\epsilon,Y_{t_0+\epsilon}^\epsilon}
 ,Y_T^{t_0+\epsilon,Y_{t_0+\epsilon}^\epsilon}\>+ o(\epsilon)
 \\ \dis   \!\!\!&=&\!\!\!
 \E
 \< P_{t_0+\epsilon}Y_{t_0+\epsilon}^\epsilon,
 Y_{t_0+\epsilon}^\epsilon\> + o(\epsilon).
\end{array}
\end{equation}
We wish to replace $P_{t_0+\epsilon}$ by $P_{t_0}$ in the above that is we claim that
 ${
  \E \< (P_{t_0+\epsilon}-P_{t_0})\, Y_{t_0+\epsilon}^\epsilon,
Y_{t_0+\epsilon}^\epsilon\>= o(\epsilon),}
$
or equivalently that ${
  \E  \< (P_{t_0+\epsilon}-P_{t_0})\, \epsilon^{-1/2}Y_{t_0+\epsilon}^\epsilon,
\epsilon^{-1/2} Y_{t_0+\epsilon}^\epsilon\>\to 0.}$

To prove the  above claim we need a compactness argument.
 A similar argument will allow us to approximate $P$ by suitable finite dimensional projections at the end of this proof.

\noindent By the Markov inequality and Proposition \ref{estimates on YZ} if we set
$K_\delta =\{ f\in L^4\;:\; f\in D(-A)^\eta,
\|f\|_{D(-A)^\eta}\le C_0\delta^{-1/4}\}$, for a suitable constant $C_0$,
and denote  by $\Omega_{\delta,\epsilon}$ the event
$\{\epsilon^{-1/2}(-A)^{-\eta} Y^\epsilon_{t_0+\epsilon}\in K_\delta\}$
we get
$
\P(\Omega_{\delta,\epsilon}^c)
\le \delta.
$
\newline\noindent Moreover
$$\begin{array}{l}
 \E  \< (P_{t_0+\epsilon}-P_{t_0})\, \epsilon^{-1/2}Y_{t_0+\epsilon}^\epsilon,
\epsilon^{-1/2} Y_{t_0+\epsilon}^\epsilon\>
\\\dis =
 \E  [\< (P_{t_0+\epsilon}-P_{t_0})\, \epsilon^{-1/2}Y_{t_0+\epsilon}^\epsilon,
\epsilon^{-1/2} Y_{t_0+\epsilon}^\epsilon\>1_{\Omega_{\delta,\epsilon}^c}]+
 \E  [\< (P_{t_0+\epsilon}-P_{t_0})\, \epsilon^{-1/2}Y_{t_0+\epsilon}^\epsilon,
\epsilon^{-1/2} Y_{t_0+\epsilon}^\epsilon\>1_{\Omega_{\delta,\epsilon}}]
\\\dis
=:A_1^\epsilon+ A_2^\epsilon.
\end{array}
$$
By the H\"older inequality
$$
|A_1^\epsilon|\le (\E\|P_{t_0+\epsilon}-P_{t_0}\|^2_\call)^{1/2}
(\E\|\epsilon^{-1/2}Y^\epsilon_{t_0+\epsilon}\|_4^8)^{1/4}\P(\Omega_{\delta,\epsilon}^c)^{1/4},
$$
and from the estimates in Proposition \ref{prop_stime-su-P}
 we conclude that
 $
|A_1^\epsilon|\le c\P(\Omega_{\delta,\epsilon}^c)^{1/4}=O(\delta^{1/4})$

\noindent On the other hand, recalling the definition of $\Omega_{\delta,\epsilon}$,
$
|A_2^\epsilon|\le
\E  \sup_{f\in K_\delta }|\< (P_{t_0+\epsilon}-P_{t_0})\, (-A)^\eta f,
(-A)^\eta f\>1_{\Omega_{\delta,\epsilon}}|.
$

\noindent
Since $K_\delta$ is compact in $L^4$, it can be covered by a finite number $N_\delta$
of open balls with radius $\delta$ and centers denoted $f_i^\delta$, $i=1,\ldots,N_\delta$.
Since $D(-A)^\eta$ is dense in $L^4$, we can assume that $f_i^\delta\in D(-A)^\eta$.
Given $f\in K_\delta$, let $i$ be such that $\|f-f_i^\delta\|_4<\delta$; then writing
$$\begin{array}{l}
  \< (P_{t_0+\epsilon}-P_{t_0}) (-A)^\eta f,(-A)^\eta f\>=
  \< (P_{t_0+\epsilon}-P_{t_0}) (-A)^\eta f_i^\delta,(-A)^\eta f_i^\delta\>
\\\dis -
  \< (P_{t_0+\epsilon}-P_{t_0}) (-A)^\eta (f-f_i^\delta),(-A)^\eta (f-f_i^\delta)\>
  +2 \< (P_{t_0+\epsilon}-P_{t_0}) (-A)^\eta f,(-A)^\eta (f-f_i^\delta)\>
\end{array}
$$
and taking expectation,
it follows from (\ref{stimaindaetaperp})
that
$$
|A_2^\epsilon|\le
\sum_{i=1}^{N_\delta}
\E |\< (P_{t_0+\epsilon}-P_{t_0}) (-A)^\eta f_i^\delta,(-A)^\eta f_i^\delta\>|
+c(T-t_0-\epsilon)^{-2\eta}
[\delta^2
  + \delta^{3/4}],
$$
and by the second statement in Proposition \ref{prop_stime-su-P} we conclude that
$$
\limsup_{\epsilon \downarrow 0} |A_2^\epsilon|\le
 c(T-t_0)^{-2\eta}
[\delta^2
  + \delta^{3/4}].
$$
Letting $\delta\to 0$ we obtain $|A_1^\epsilon|+|A_2^\epsilon|\to 0$ and the proof
that  ${
  \E  \< (P_{t_0+\epsilon}-P_{t_0})\, \epsilon^{-1/2}Y_{t_0+\epsilon}^\epsilon,
\epsilon^{-1/2} Y_{t_0+\epsilon}^\epsilon\>\to 0}$
is finished.

We come now to show that
$${
  \E  \<  P_{t_0}  Y_{t_0+\epsilon}^\epsilon,
Y_{t_0+\epsilon}^\epsilon\>= \E\int_{t_0}^{t_0+\epsilon}
 \< P_s \delta^\epsilon \sigma_j(s,\cdot),\delta^\epsilon \sigma_j(s,\cdot)\>\,ds
+ o(\epsilon).}
$$
 If  we treat $A$ and $P_{t_0}$ as bounded operators in $H$ we  get,  by It\^o rule:
$$ \E  \<  P_{t_0}  Y_{t_0+\epsilon}^\epsilon,
Y_{t_0+\epsilon}^\epsilon\>=2\E \int_{t_0}^{t_0+\epsilon}\!\!\!\!\<  P_{t_0}  Y_{s}^\epsilon,
(A+ b'(\bar{X}_t,\bar{u}_t)) Y_{s}^\epsilon\> ds+ \E\int_{t_0}^{t_0+\epsilon}\!\!\!\!
 \< P_{t_0} \delta^\epsilon \sigma_j(s,\cdot),\delta^\epsilon \sigma_j(s,\cdot)\>\,ds$$
and the claim follows recalling that  $\E  |Y^\epsilon_t|^{2}=O( \epsilon)$ and   the {``continuity"} of $P$
stated in Proposition \ref{prop_stime-su-P}.

\noindent The  general case is more technical and requires a double approximation:
 $A$ by its Yosida Approximations and
 $P$ by finite dimensional projections  $
 { P_t^N(\omega)f:
=\sum_{i,j=1}^N \<P_t(\omega)e_i,e_j\> \<e_i,f\>_2 e_j,}$
$ f\in L^4,$
where  $(e_i)_{i\ge 1}$ is an  orthonormal basis in
$L^2$ which is also a Schauder basis of $L^4$.

\noindent The conclusion of the proof of the maximum principle is now standard (see, e.g. \cite{HuPe}, \cite{TaLi} or \cite{Peng}). We just have to write
 $J(u^\epsilon)-J(u)$ using (\ref{penultima_stima}), (\ref{delta_3}) and (\ref{eq-penultima}), to recall that
 $0\leq \epsilon^{-1}( J(u^\epsilon)-J(u))$ and to let
$\epsilon\rightarrow 0$.

\end{document}